\documentclass[a4,12pt]{article}

\usepackage{amsmath,amssymb}

\newtheorem{thm}{Theorem}
\newtheorem{prop}{Proposition}
\newtheorem{lem}{Lemma}
\newtheorem{cor}{Corollary}

\title{A characterization of \\ Gorenstein toric Fano $n$-folds with index $n$ \\
 and Fujita's conjecture}
\author{Shoetsu Ogata\thanks{e-mail:  ogata{\char'100}math.tohoku.ac.jp} and  
Huai-Liang Zhao\\
Mathematical Institute, Tohoku University\\ Sendai 980-8578, Japan}

\begin{document}
\maketitle

\begin{abstract}
We give a characterization of Gorenstein toric Fano varieties of dimension $n$
with index $n$  among toric varieties.  As an application, we  give a stronger version of 
Fujita's freeness conjecture  and  also give a simple proof of
Fujita's very ampleness conjecture on Gorenstein toric varieties.
\end{abstract}
\section*{Introduction}

A nonsingular projective variety $X$ is called {\it Fano} if its anti-canonical divisor
$-K_X$ is ample.  For a Fano variety $X$ 
  the number $i_X:=\max\{i\in \mathbb{N}; -K_X=iD 
\ \mbox{for ample}
\ \mbox{ divisor $D$}\}$ is called the {\it Fano index}, or simply {\it index}.
Kobayashi and Ochiai \cite{KO}
showed that  a Fano variety $X$ of dimension $n$ with index $n+1$
 is the projective
$n$-space and that a Fano variety with index  $n$ is the hyperquadric.

First, we give a characterization of the projective space among toric varieties.

\begin{thm}\label{0;t0}
Let $X$ be a projective toric variety of dimension $n$ $(n\ge2)$.
We assume that there exists an ample line bundle $L$ on $X$ with $\dim \Gamma(X,
L^{\otimes n}\otimes \omega_X)=0$,
where $\omega_X$ is the dualizing sheaf of $X$.
Then $X$ is the projective space of dimension $n$.
\end{thm}
This is not a new result.  If we consider the lattice polytope corresponding to a polarized
toric variety $(X, L)$, this is a characterization of the basic lattice simplex 
treated by Batyrev and Nill \cite[Proposition 1.4]{BN}.
In this article, we use the normality of some multiple of lattice polytope in order to
characterize Gorenstein toric Fano varieties.   In our purpose, this characterization stated as
Lemma~\ref{l3} in Section~\ref{sect1} is essencial. 
Theorem~\ref{0;t0} is obtained as Corollary~\ref{cor1} in Section~\ref{sect1}.

We also give a characterization of 
 a Gorenstein toric Fano variety of dimension $n$ with index $n$
 among toric varieties.

\begin{thm}\label{0;t1}
Let $X$ be a projective toric variety of dimension $n$ $(n\ge2)$.
We assume that there exists an ample line bundle $L$ on $X$ with $\dim \Gamma(X,
L^{\otimes n}\otimes \omega_X)=1$.
Then $X$ is a Gorenstein toric Fano variety of dimension $n$ with  index $n$.
\end{thm}

Batyrev and Nill obtain a classification of lattice polytope $P$ of dimension $n$ such that
the multiple $(n-1)P$ does not contain lattice points in its interior \cite[Theorem 2.5]{BN}.
This can be interpreted as a classification of polarized toric variety $(X, L)$ with
$\Gamma(X,L^{\otimes (n-1)}\otimes \omega_X)=0$.
We do not use their classification but do an elementary argument about the shape
of polytopes. 

In this paper, we call  a variety of dimension $n$, simply, an $n$-fold.
We see that a Gorenstein toric Fano surface with index 2 is the quadratic surface
$\mathbb{P}^1\times\mathbb{P}^1$ or the weighted projective plane 
$\mathbb{P}(1,1,2)$.
For higher dimension $n\ge3$, we see that the weighted projective $n$-space
$\mathbb{P}(1,1,2,\dots, 2)$ is a Gorenstein toric Fano $n$-fold with  index $n$.
We also have another toric Fano $n$-folds with  index $n$.

\begin{thm}\label{0;t2}
Let $X$ be a Gorenstein toric Fano variety of dimension $n$ $(n\ge2)$ with index $n$.
Let $D$ be an ample Cartier divisor on $X$ with $-K_X=nD$.
Then $D$ is very ample and $(X, D)$ is a quadratic hypersurface in $\mathbb{P}^{n+1}$
which is a cone over a plane conic, i.e., the weighted projective space $\mathbb{P}(1,1,2,\dots, 2)$, or 
a cone over the quadratic surface  $\mathbb{P}^1\times\mathbb{P}^1\hookrightarrow
\mathbb{P}^3$.
\end{thm}

We may expect the same characterization of a Gorenstein toric Fano $n$-fold with
index $n-1$ among toric varieties for $n\ge3$.  We give a special case as Proposition~\ref{pr1}
in Section~\ref{sect5}.

\begin{thm}\label{0;t3}
Let $X$ be a projective toric variety of dimension $n$ $(n\ge3)$.
We assume that there exists an ample line bundle $L$ on $X$ with 
$\dim \Gamma(X, L)=n+1$ and 
$\dim \Gamma(X,L^{\otimes (n-1)}\otimes \omega_X)=1$.
Then $X$ is a Gorenstein toric Fano variety of dimension $n$ with  index $n-1$.
\end{thm}

Batyrev and Juny \cite{BJ} classified Gorenstein toric Fano $n$-folds with index $n-1$.
For the proof of Theorem~\ref{0;t3} we do not use the classification of Batyrev and Juny.
We also remark that there exists a non-Gorenstein toric $n$-fold $X$ with the ample
line bundle $L$ such that $\dim \Gamma(X, L)=n+1$ and 
$\dim \Gamma(X,L^{\otimes (n-2)}\otimes \omega_X)=1$ for every $n\ge4$.
We give  examples in the last of this article.

As a corollary of above theorems, we obtain a strong version of Fujino's theorem\cite{Fn}
for Gorenstein toric varieties.  This is given in Section~\ref{sect4}.
Theorem of this type is called "Fujita's freeness conjecture" \cite{Fj} in general.

\begin{thm}\label{0;tm1}
Let $X$ be a Gorenstein projective toric variety of dimension $n$ $(n\ge2)$.
We assume that  $X$ is not the projective space, a $\mathbb{P}^{n-1}$-bundle
over $\mathbb{P}^1$ nor a Gorenstein toric Fano $n$-fold with index $n$.
If an ample line bundle $L$ on $X$ satisfies
that the intersection number with every irreducible invariant curve is at least $n-1$,
then the adjoint bundle $L+K_X$ is nef.
\end{thm}

We also give a simple proof of the following theorem, which is a special case of
theorem of Payne \cite{P} for Gorenstein toric varieties.  
This is given in Section~\ref{sect5}.
Theorem of this type is also called "Fujita's very ampleness conjecture" in general.

\begin{thm}\label{0;tm2}
Let $X$ be a Gorenstein projective toric variety of dimension $n$ not isomorphic to
the projective space.
If an ample line bundle $L$ on $X$ satisfies
that the intersection number with every irreducible invariant curve is at least $n+1$,
then the adjoint bundle $L+K_X$ is very ample.
\end{thm}

We note that the condition on the intersection number "at least $n+1$"
is trivially best possible for Gorenstein toric Fano $n$-folds with index $n$.
We remark, however, that even if we make an exception on $X$ as in  Theorem~\ref{0;tm1},
we cannot weaken the condition  "at least $n+1$".
  We give  examples in Section~\ref{sect5} for all dimension $n\ge3$.

\section{Toric Varieties and Lattice Polytopes}\label{sect1}

In this section, we recall the correspondence between polarized toric varieties and 
lattice polytopes and give criterions for  polytopes to be basic.

Let $M=\mathbb{Z}^n$ be a free abelian group of rank $n$ and $M_{\mathbb{R}}:=
M\otimes_{\mathbb{Z}}\mathbb{R} \cong\mathbb{R}^n$ the extension of coefficients
into real numbers.
We define a {\it lattice polytope} $P$  in $M_{\mathbb{R}}$ as the convex hull 
$P:= \mbox{Conv}\{m_1,
\dots, m_r\}$ of a finite subset $\{m_1, \dots, m_r\}$ of $M$.
We define the dimension of a lattice polytope $P$ as that of the smallest
affine subspace containing $P$.
We sometimes call an $n$-polytope a polytope of dimension $n$.

The space $\Gamma(X, L)$ of global sections of an 
ample line bundle $L$ on a toric variety $X$ of dimension $n$
 is parametrized by the set of lattice points in a lattice polytope $P$
of dimension $n$ (see, for instance, Oda's book\cite[Section 2.2]{Od} 
or Fulton's book\cite[Section 3.5]{Fu}). 
 And $k$ times tensor product $L^{\otimes k}$ corresponds to the polytope
 $kP:=\{kx\in M_{\mathbb{R}}; x\in P\}$.   Furthermore, the surjectivity of
 the multiplication map $\Gamma(X, L^{\otimes k})\otimes \Gamma(X, L) \to
 \Gamma(X, L^{\otimes(k+1)})$
 is equivalent to the equality
 \begin{equation}\label{eq1}
 (kP)\cap M +P\cap M = ((k+1)P)\cap M.
 \end{equation}
In particular, the dimension of the space of global sections of $L\otimes \omega_X$ is
equal to the number of lattice points contained in the interior of $P$, i.e.,
$$
\dim \Gamma(X, L\otimes \omega_X)=\sharp \{(\mbox{Int}P)\cap M\}.
$$
A lattice polytope $P$ is called {\it normal} if the equality~(\ref{eq1})
holds for all $k\ge1$.

By using the generalized Castelnuovo's Lemma of Mumford \cite{Mf},
we can prove the following lemma.
We can also prove this in terms of lattice polytopes.
\begin{lem}\label{l1}
Let $P$ be a lattice polytope of dimension $n$.
If there exists an integer $r$ with $1\le r\le n-1$ satisfying the condition that
the multiple $rP$ does not contain lattice points in its interior,
then the equality
$$
(kP)\cap M +P\cap M =((k+1)P)\cap M
$$
holds for all integers $k\ge n-r$.
\end{lem}

Even if $P$ contains lattice points in its interior,
we have a   result of Nakagawa \cite{N}.

\begin{thm}[Nakagawa\cite{N}]\label{2;t1}
If a lattice polytope $P$ is dimension $n$, then  the equality
$$
(kP)\cap M +P\cap M =((k+1)P)\cap M
$$
holds for all integers $k\ge n-1$.
\end{thm}

Let $\{e_1, \dots, e_n\}$ be a $\mathbb{Z}$-basis of $M$. 
A lattice polytope is called {\it basic} if it is isomorphic to the $n$-simplex
$$
\Delta_0^n:=\mbox{Conv}\{0, e_1, \dots, e_n\}
$$
by a unimodular transformation of $M$.
For a basic $n$-simplex $P$ we see that $\mbox{Int}(nP)\cap M=\emptyset$ and
$\sharp \mbox{Int}((n+1)P)\cap M=1$.
A lattice polytope $P$ of dimension $n$ is called an {\it empty  lattice $n$-simplex}
if the set of lattice points in $P$ are only $n+1$ vertices.

\begin{lem}\label{l2}
Let $P$ be an empty lattice  $n$-simplex.
If the multiple $(n-1)P$ does not contain lattice points in its interior,
then it is basic.
\end{lem}
{\it Proof}.
From Lemma~\ref{l1}, we see that $P$ is normal, that is, the equality~(\ref{eq1})
holds for all $k\ge1$.

By taking a coordinates of $M$ so that the origin is a vertex of $P$,
set $P=\mbox{Conv}\{0, m_1, \dots, m_n\}$.
Set the cone 
$$
C(P):=\mathbb{R}_{\ge0}m_1+\dots +\mathbb{R}_{\ge0}m_n \subset M_{\mathbb{R}}.
$$
Then the normality of $P$ implies that the semi-group $C(P)\cap M$ is generated
by $m_1, \dots, m_n$.   This proves the Lemma. \hfill $\Box$

\medskip

\begin{lem}\label{l3}
Let $P$ be a lattice $n$-polytope.
If the multiple $nP$ does not contain lattice points in its interior,
then it is basic.
\end{lem}
{\it Proof}.   Take a coordinates of $M$ so that the origin coincides with a vertex of $P$.

First we treat the case that $P$ is a lattice $n$-simplex $\mbox{Conv}\{0,m_1, \dots,
m_n\}$.
If $P$ is an empty lattice simplex, then it is basic by Lemma~\ref{l2}.
We assume that a face $F=\mbox{Conv}\{0, m_1, \dots, m_r\}$ has a lattice point $m$
 in its relative interior.  Then $m+m_{r+1}+\dots +m_n$ is contained in the interior of $nP$.
This contradicts to the assumption.  If a face $F'=\mbox{Conv}\{m_1, \dots, m_r\}$ 
contains an lattice point $m$ in its interior,
then $m+m_{r+1}+\dots +m_n$ is contained in the interior of $nP$.
Thus we see that $P$ is an empty lattice simplex if it is an lattice simplex.

Next we treat the general case.
Let $C_0(P):=\mathbb{R}_{\ge0}(P)$ be the cone of $P$ with apex $0$.
If a lattice point in $P$ is contained in the interior of the cone $C_0(P)$, then
it is contained in the interior of $2P$. 

When $n=2$, we see $P$ is  a triangle,  hence, it is basic by the first part of the proof.

In general dimension $n\ge3$, let $\{m_1, \dots, m_t\}$ be the set of the end points of all edges 
of $P$ through the origin.  
We may assume that $\{m_1, \dots, m_n\}$ is linearly independent.
Set $Q:=\mbox{Conv}\{0, m_1, \dots, m_n\}$.
Then $Q$ is basic by the first part because $\dim Q=n$ and 
$\mbox{Int}(nQ)\cap M=\emptyset$.

If $t>n$, that is, $P\not=Q$,
then we have a facet, say, $F=\mbox{Conv}\{0, m_1, \dots,$   $m_{n-1}\}$
of $Q$ which is not a facet of $P$.   Since 
$\mbox{Int}(n)F\cap M$ contains a lattice point, it would be a  lattice point  in the interior
of $nP$.  This contradicts to the assumption.  
Hence, we see that $t=n$, that is, $P$ is a lattice $n$-simplex.
By the first part of this proof, it is basic.  \hfill $\Box$

\medskip
We may consider Lemma~\ref{l3} as a characterization of the projective $n$-space
among toric varieties.

\begin{cor}\label{cor1}
Let $X$ be a projective toric variety of dimension $n$ $(n\ge2)$.
We assume that there exists an ample line bundle $L$ on $X$ with $\dim \Gamma(X,
L^{\otimes n}\otimes \omega_X)=0$.
Then $X$ is the projective space $\mathbb{P}^n$ and $L=\mathcal{O}(1)$.
\end{cor}

\section{Gorenstein Fano with index $n$}\label{sect2}

In this section, we will give a characterization of a Gorenstein toric Fano $n$-fold
with index $n$ among toric varieties.

Let $M$ be a free abelian group of rank $n\ge2$ and $P\subset M_{\mathbb{R}}$ 
a lattice $n$-polytope.
For a vertex $v$ of $P$, make the cone 
$$
C_v(P):=\mathbb{R}_{\ge0}(P-v)=\{r(x-v)\in M_{\mathbb{R}}; r\ge0\ \mbox{and
$x\in P$}\}.
$$
We call $P$ {\it Gorenstein at} $v$ if there exists a lattice point $m_0$ in $C_v(P)$
such that the equality
$$
(\mbox{Int} C_v(P))\cap M= m_0+C_v(P)\cap M
$$
holds.  We call $P$ {\it Gorenstein} if it is Gorenstein at all vertices.

We also define the notion of very ampleness.
We call $P$  {\it very ample} at $v$ if the semi-group $C_v(P)\cap M$
is generated by $(P-v)\cap M$ and $P$ {\it very ample} if it is very ample at all vertices.

Let $\{e_1, \dots, e_n\}$ be a $\mathbb{Z}$-basis of $M$.
Set
\begin{eqnarray*}
P_n:&=&\mbox{Conv}\{0, 2e_1, e_2, \dots, e_n\},\\
Q_n:&=& \mbox{Conv}\{0, e_1, e_2, e_1+e_2, e_3, \dots, e_n\}.
\end{eqnarray*}
Then they are a very ample Gorenstein polytopes and $\sharp\{\mbox{Int}(nP_n)\cap M \}=
\sharp\{\mbox{Int}(nQ_n)\cap M\}=1$.
The polytope $P_n$ corresponds to the polarized variety $(\mathbb{P}(1,1,2,\dots,2), 
\mathcal{O}(2))$.
The polygon $Q_2$ corresponds to $(\mathbb{P}^1\times \mathbb{P}^1, \mathcal{O}(1,1))$,
and $Q_3$ corresponds to a cone over the quadratic surface $\mathbb{P}^1\times \mathbb{P}^1
\hookrightarrow \mathbb{P}^3$.
\medskip

\begin{prop}\label{p1}
Let $P\subset M_{\mathbb{R}}$ be a lattice $n$-polytope for $n\ge2$
satisfying the condition $\sharp\{\mbox{\rm Int}(nP)\cap M\}=1$.
Then $P$ is isomorphic to $P_n$ or $Q_n$ by a unimodular transformation of $M$.
\end{prop}

For the proof we fix notations.  Take a vertex $v\in P$.  We may set 
 so that the origin coincides with  $v$ by a parallel transformation of $M$.
Let $\{m_1, \dots, m_t\}$ be the set of the end points of all edges through the origin
and define as
\begin{equation}\label{2;e1}
Q:=\mbox{Conv}\{0, m_1, \dots, m_t\}.
\end{equation}
Set $\{\tilde{m}\}=\mbox{Int}(nP)\cap M$.
We note that $(n-1)P$ does not contain lattice points in its interior
because if $m'\in M$ is in the interior of $(n-1)P$, then $\{m'+m_i;\ i=1, \dots, t\}$
is in the interior of $nP$.
Hence $\mbox{Int}(n-1)Q\cap M=\emptyset$.

First consider the case when $n=2$.  Then $t=n=2$.
Since $\mbox{Int}P\cap M=\emptyset$,
we may set $m_1=ae_1, m_2=be_2$ with positive integers $a\ge b$.
When $P$ is a triangle, we see that $a=b=2$ or $b=1$, and the condition
$\sharp\{\mbox{Int}(2P)\cap M\}=1$ implies $a=2, b=1$.
Thus $P\cong P_2$.

If $P$ is not a triangle, then $b=1$ and $P$ is a quadrangle with the other vertex
$ce_1+e_2$ for $c\ge1$.   Since $2P$ contains only one lattice point in its interior,
we see $a=c=1$.  This $P$ coincides with $Q_2$.

In the following we set $n\ge3$.

\begin{lem}\label{2;l1}
Let $Q$ be the $n$-polytope defined in {\rm (\ref{2;e1})}.
Assume that $\mbox{\rm Int}(nQ)\cap M=\{\tilde{m}\}$.
Then $Q$ is isomorphic to $P_n$ or $Q_n$.
\end{lem}
{\it Proof}.
Since $\tilde{m}$ is contained in the interior of $C_v(P)=C_v(Q)$, 
by the Carath\'eodory's theorem
we may choose $n$ elements from $m_1, \dots, m_t$ such that the simplicial cone
spanned by them contains $\tilde{m}$ and furthermore we can choose
$m_1, \dots, m_r$ so that the subcone spanned by them contains $\tilde{m}$ in
its relative interior by renumbering.

Set $G:=\mbox{Conv}\{m_1, \dots, m_r\}$ and $\Tilde{G}:=\mbox{Conv}\{0, 
{G}\}$. Then $\dim\Tilde{G}=r$.  Since $\mbox{Int}(n-1)\Tilde{G}\cap M=\emptyset$
and $\mbox{Int}(r+1)\Tilde{G}\cap M\not=\emptyset$, we have $n-1<r+1$,
hence, $r\ge n-1$.
We separate our argument into two cases.

The case (a): $r=n-1$.  
Since $\mbox{Int}(n-1)\Tilde{G}\cap M=\emptyset$, we see that $\Tilde{G}$ is
basic from Lemma~\ref{l3}.
Let $H$ be the hyperplane in $M_{\mathbb{R}}$ containing $\Tilde{G}$.
We can divide $Q$ into the union of two $n$-polytopes $Q_{(1)}$ and 
$Q_{(2)}$ separated by $H$.
Since $\tilde{m}\in n\tilde{G}$ and $\sharp\{\mbox{\rm Int}(nQ)\cap M\}=1$, 
we see that $\mbox{Int}(nQ_{(i)})\cap M=\emptyset$ for $i=1, 2$.
Thus $Q_{(i)}$ is basic.  We may set $m_j=e_j$ for $j=1, \dots, n$ so that $Q_{(1)}=
\mbox{Conv}\{\tilde{G}, e_n\}$, $Q_{(2)}=\mbox{Conv}\{\tilde{G}, m_{n+1}\}$
and $t=n+1$.
Furthermore, we set 
$$
m_{n+1}=a_1e_1+\dots +a_{n-1}e_{n-1} -e_n.
$$
Since $Q$ is convex, we see that $a_i\ge0$ and $1\le\sum_i a_i\le2$.
If $a_2=\dots=a_{n-1}=0$, then $m_1=e_1$ would not be a vertex of $Q$.
Hence, we have $a_1=a_2=1, a_3=\dots=a_{n-1}=0$.
If we set $e'_i=e_i-m_{n+1}$ for $i=1, \dots, n-1$ and $e'_n=0-m_{n+1}$, then 
$e_n-m_{n+1}=e'_1+e'_2$ and $\{e'_1, \dots, e'_n\}$ is a $\mathbb{Z}$-basis
of $M$.  Thus, we see that  $Q$ is isomorphic to $Q_n$.

The case (b): $r=n$.
If $t>n$, then there is a facet, say, 
$F:=\mbox{Conv}\{0, m_1, \dots, m_{n-1}\}$ of $\Tilde{G}$ which 
is not facet of $Q$, that is, the relative interior of $F$ is contained in the interior of $Q$.  
Since $\mbox{Int}(nF)\cap M\not=\emptyset$, the interior of $nQ$ would contain 
lattice points more than one.  Thus we have $t=n$.

Take the primitive elements $m'_i\in M$ for $i=1, \dots, n$ 
so that $m_i=a_im'_i$ with positive integers $a_i$.
Set $Q':=\mbox{Conv}\{0, m'_1, \dots, m'_n\}$.
We may set $a_1\ge  a_2 \ge\dots \ge a_n\ge1$. For $2\le r\le n-1$, set
$$
F_r:=\mbox{Conv}\{0, m_1,  \dots, m_r\}.
$$
If $F_2$ contains lattice points in its relative interior, then the interior of
$(n-1)Q$ would contain lattice points.  Thus $a_1=a_2=2$, or $a_2=1$.
If $a_1=a_2=2$, then $2F_2$ contains lattice points
more than two in its relative interior,
hence, the interior of $nQ$ would contain lattice points more than two.  Thus we have
$a_2=1$.

When $a_1\ge2$, set
$$
G:=\mbox{Conv}\{m'_1, m_2, \dots, m_n\}.
$$
Since the relative interior of $G$ is contained in the interior of $Q$ and $\dim G=n-1$,
the interior of 
$(n-1)G$ does not contain lattice points, hence, $G$ is basic from Lemma~\ref{l3}.
Since the number of lattice points in the interior of $nG$ is one, we see $a_1=2$.
Set
$$
\Tilde{G}:=\mbox{Conv}\{0, m'_1, m_2, \dots, m_n\}.
$$
Since $n\Tilde{G}$ could not contain lattice points in its interior, it is basic.
Thus $Q$ is isomorphic to $P_n$.

When $a_1=1$, we see that $Q$ is not basic because of the assumption
$\mbox{Int}(nQ)\cap M\not=\emptyset$.
If $Q$ is an empty lattice simplex, then the condition $\mbox{Int}(n-1)Q\cap M=\emptyset$
implies that $Q$ is basic by Lemma~\ref{l2}.
Thus there exists a lattice point $v$ in $Q$ other than its vertices.
After renumbering, assume that  $F_r$ contains $v$ with the smallest dimension.
We note $F_r$ does not contain lattice points in its interior.
Set
$$
F'_r:=\mbox{Conv}\{m_1,  \dots, m_r\}.
$$
Since $v$ is contained in the relative interior of $F'_r$, 
it is contained in the relative interior of $2F_r$, hence it is 
contained in that of $(2+n-r)Q$.  Thus we have $r=2$.
In this case, we see that $Q$ is isomorphic to $P_n$.
\hfill $\Box$

\medskip
\noindent
{\it Proof of Proposition~\ref{p1}}.
We separate our argument into two cases according to the existence of interior lattice
points of $nQ$.

The case I: $\mbox{Int}(nQ)\cap M=\emptyset$.
Then $Q$ is basic by Lemma~\ref{l3} and $P\not=Q$.
Let $v'$ be a vertex of $P$ not contained in $Q$.  If $v'$ is contained in the interior of
the cone $C_v(P)$, then it would be contained in the interior of $2P$.
Thus $v'$ is contained in the relative interior of a face cone
$$
\mathbb{R}_{\ge0}m_1+\dots +\mathbb{R}_{\ge0}m_r
$$
with $2\le r\le n-1$ by renumbering if necessary.
Let $G\subset P$ be the face containing $m_1, \dots , m_r$ and $v'$.
Then $v'$ is contained in the relative interior of $2G$ and $v'+m_{r+1}+\dots +m_n$
is also contained in the interior of $(n-r+2)P$, hence, $r\le2$ and $r=2$.
Since $(\mbox{\rm Int}(n-1)P)\cap M=\emptyset$, the face $G$ does not contain
lattice points in its relative interior.  Since $\{m_1, \dots, m_n\}$ is a $\mathbb{Z}$-basis
of $M$, we can write as 
$$
v'=am_1+m_2
$$
with a positive integer $a$.
If $a\ge2$, then the number of lattice points in the interior of $nP$ would be more than one.
Thus $a=1$ and $P$ is isomorphic to $Q_n$.

The case II: $\mbox{Int}(nQ)\cap M\not=\emptyset$.  Since $\sharp\{\mbox{\rm Int}(nP)\cap M\}=1$, we have $\mbox{Int}(nQ)\cap M=\{\tilde{m}\}$.
If $P\not=Q$, then there exists a facet $F$ of $Q$ which is not a facet of $P$.
Since $\mbox{Int}(nF)\cap M\not=\emptyset$, they would be contained in the interior of $nP$
and does not coincide with $\tilde{m}$,
which contradicts to the assumption.  Thus we have
 $P=Q$.  This case has been proved in Lemma~\ref{2;l1}.
\hfill $\Box$

\medskip
We note that Proposition~\ref{p1} is an interpretation of Theorems~\ref{0;t1} and \ref{0;t2}
in terms of lattice polytopes.

\section{A certain class of Gorenstein $n$-polytopes}\label{sect3}

In this section, we  determine an $n$-polytope $P$ such that $(n-1)P$ does not
contain lattice points in its interior. 
This polytope will play an important role in the next section.

We recall that a basic $n$-simplex is isomorphic to
$$
\Delta_0^n=\mbox{Conv}\{0, e_1, \dots, e_n\}.
$$
Now for positive integers $a_1\le a_2\le \dots \le a_n$ set
$$
R_n:=\mbox{Conv}\{\Delta_0^{n-1}, e_1+a_1e_n, \dots, e_{n-1}+a_{n-1}e_n, a_ne_n\}.
$$
Then $R_n$ is a nonsingular $n$-polytope and $(\mbox{Int}(n-1)R_n)\cap M=\emptyset$.
We also know $(\mbox{Int}(n-1)P_n)\cap M=(\mbox{Int}(n-1)Q_n)\cap M=\emptyset$.

\medskip

\begin{prop}\label{p1a}
Let $P\subset M_{\mathbb{R}}$ be a Gorenstein $n$-polytope for $n\ge3$
satisfying the condition $(\mbox{\rm Int}(n-1)P)\cap M=\emptyset$
and $(\mbox{\rm Int}(nP))\cap M\not=\emptyset$.
Then $P$ is isomorphic to $P_n$, $Q_n$ or $R_n$ by a unimodular transformation of $M$.
\end{prop}

As in the proof of Proposition~\ref{p1}, we take a vertex $v$ of $P$
so that $v$ coincides with the origin $0$ of $M$.
Since $P$ is Gorenstein, the cone $C_v(P)=\mathbb{R}_{\ge0}(P-v)$ contains the
lattice point $m_0$ in its interior such that
\begin{equation}\label{3;e1}
(\mbox{Int}C_v(P))\cap M = m_0+ C_v(P)\cap M.
\end{equation}
Let $\{m_1, \dots, m_t\}$ be the set of the end points of all edges through 
the vertex $v=0$
and $Q:=\mbox{Conv}\{0, m_1, \dots, m_t\}$.

We give a lemma useful in the following sections.

\begin{lem}\label{3;l0}
Let $Q$ be the $n$-polytope defined above.
Assume that 
$\mbox{\rm Int}(rQ)\cap M\not=\emptyset$ for some $1\le r\le n$.
Then $m_0$ is contained in the interior of $rQ$.
\end{lem}
{\it Proof}.
Since $m_0$ is contained in the interior of the cone $C_v(Q)$, by the Carath\'eodory's Theorem
we can choose $n$ elements from $m_1, \dots, m_t$ such that the simplicial cone
spanned by them contains $m_0$.  By renumbering, we assume that the cone
$C(m_1, \dots, m_n)$ spanned by $m_1, \dots, m_n$ contains $m_0$.
Let $\nu\in M^{\star}_{\mathbb{Q}}$ be the rational point in the dual space to
$M_{\mathbb{Q}}$ such that the hyperplane $H_1:=\{x\in M_{\mathbb{R}};
\langle \nu, x\rangle =1\}$ contains all $m_1, \dots, m_n$.
We note that the parallel hyperplane $H_0=\{\langle \nu, x\rangle=0\}$
touches the cone $C_v(Q)$ at one point $0$, in particular, we have
$\langle \nu, x\rangle\ge0$ for all $x\in C_v(Q)$.

Set $\tilde{m}\in \mbox{Int}(rQ)\cap M$.  Then $\langle \nu, \tilde{m}\rangle <r$.
If $m_0$ is not contained in the interior of $rQ$. then $\langle \nu, m_0\rangle\ge r$.
Since $C_v(Q)$ is Gorenstein at $0$, there exists $m'\in (nQ)\cap M$ such that
$\tilde{m}=m_0+m'$ by the equation~(\ref{3;e1}).  
But $\langle \nu, m'\rangle <0$.  This is a contradiction.
Thus we have $m_0\in \mbox{Int}(rQ)$.  \hfill $\Box$

\medskip

We also define a Gorenstein Fano $n$-polytope $Q'_n$ isomorphic to $Q_n$ as
$$
Q'_n:=\mbox{Conv}\{0, e_1, \dots, e_n, e_1+e_2-e_{n}\}.
$$

\begin{lem}\label{3;l1}
Let $Q$ be the $n$-polytope defined above.
Assume that $\mbox{\rm Int}(n-1)Q\cap M=\emptyset$
and that $\mbox{\rm Int}(nQ)\cap M\not=\emptyset$.
Then $Q$ is isomorphic to $P_n$ or $Q'_n$.
\end{lem}
{\it Proof}.
Since $m_0$ is contained in the interior of $C_v(P)$, as in the proof of Lemma~\ref{2;l1}
we  can choose
$m_1, \dots, m_r$ so that the  simplicial cone spanned by them contains $m_0$ in
its relative interior by renumbering.

Set $G:=\mbox{Conv}\{m_1, \dots, m_r\}$ and $\Tilde{G}:=\mbox{Conv}\{0, 
{G}\}$.  
Then we have $r= n-1$ or $r=n$.
We separate our argument into two cases.

The case (a): $r=n-1$.  
Since $\mbox{Int}(n-1)\Tilde{G}\cap M=\emptyset$, we see that $\Tilde{G}$ is
basic from Lemma~\ref{l3}.  We note that $m_0=m_1+\dots+m_{n-1} \in
\mbox{Int}(n\tilde{G}) \subset \mbox{Int}(nQ)$.
As in the proof of Lemma~\ref{2;l1}, take
 the hyperplane $H$ in $M_{\mathbb{R}}$ containing $\Tilde{G}$, and
 divide $Q$ into the union of two $n$-polytopes $Q_{(1)}$ and 
$Q_{(2)}$ separeted by $H$.

If $\mbox{Int}(nQ_{(i)})\cap M=\emptyset$ for $i=1, 2$, then $Q=Q'_n\cong Q_n$
as in the proof of Lemma~\ref{2;l1}.

Assume that $\mbox{Int}(nQ_{(1)})\cap M\not=\emptyset$.
Set $\tilde{m}\in \mbox{Int}(nQ_{(1)})\cap M$.
Since $Q$ is Gorenstein at $0$ and since $m_0$ is contained in the boundary of
$(n-1)\tilde{G}$, hence, in that of $(n-1)Q_{(1)}$,
there exists $m'\in (Q_{(1)}\setminus H)\cap M$ with
$\tilde{m}=m_0+m'$.
Consider the $n$-polytope $R=\mbox{Conv}\{m', \tilde{G}\}$.
We see that $\tilde{m}$ is not contained in the interior of $nR$.
Since $\tilde{m}$ is in the interior of $nQ_{(1)}$, the lattice point $m'$ must
be contained in the interior of $Q_{(1)}$.  This is a contradiction.
Thus we see that $Q_{(i)}$ are basic.

The case (b): $r=n$.  From Lemma~\ref{3;l0}, we have $m_0\in\mbox{Int}(n\tilde{G})$.
We note that $\mbox{Int}(n-1)\tilde{G}\cap M=\emptyset$ by assumption.

First, consider the case that $t=n$, that is, $Q=\tilde{G}$.
We claim that 
$\mbox{Int}(n\tilde{G})\cap M=\{m_0\}$. 
If $\tilde{m}\in M$ is a lattice point in the interior of $n\tilde{G}$ other than
$m_0$, then
 there would exist a lattice point $m'\in \tilde{G}$ with $\tilde{m}=m_0+m'$ and $m'$ 
 would be  in the interior of $\tilde{G}$.   
 This contradicts to the assumption.  We confirm the claim.
 In this case, we know that $Q=\tilde{G}=P_n$ from Lemma~\ref{2;l1}.

Next, consider the case $t>n$.  
Then there is a facet, say, 
$F_{n-1}:=\mbox{Conv}\{0, m_1, \dots, m_{n-1}\}$ of $\tilde{G}$ which 
is not a facet of $Q$ and  is basic because of 
$\mbox{Int}(n-1)F_{n-1}\cap M=\emptyset$.
We can choose a new $\mathbb{Z}$-basis $\{e'_1, \dots, e'_n\}$ of $M$ as
$m_i=e'_i$ for $i=1, \dots, n-1$ and $m_n=a_1e'_1+\dots +a_ne'_n$
with $a_i\ge0$.   We note that $a_n\ge2$ because $\tilde{G}$ is not basic.
 Since $\mbox{Int}(n-1)\tilde{G}\cap M=\emptyset$, we see that $\tilde{G}$
 is normal from Lemma~\ref{l1}.  Thus $\tilde{G}$ contains a lattice point $u$
 whose $n$-th coordinate is one.
 
 Let $\tilde{F'_r}:=\mbox{Conv}\{0, m_1, \dots, m_{r-1}, m_n\}$ be the face
 containing $u$ with the smallest dimension.  If $u$ is in the relative interior, then
 we see $r=1$ by the assumption.  Then $m_n=a_ne'_n$ and $u=e'_n$.
In this case, $\tilde{G}$ is nonsingular at $0$ and $m_0=m_1+\dots+m_{n-1}+e'_n
\in \mbox{Int}(n\tilde{G})\cap M$.
The point $m_1+\dots+m_{n-1}$ is contained in the interior of $nQ$.  The equation~(\ref{3;e1})
implies that $-e'_n\in Q$,
which contradicts to that $0$ is a vertex of $Q$.

If $u$ is in the relative interior of the face
$F'_r:=\mbox{Conv}\{ m_1, \dots, m_{r-1}, m_n\}$ of $\tilde{F'_r}$,
then $r=2$ and $m_n=e'_1+a_ne'_n$, $u=e'_1+e'_n$.
When $a_n=2$,  $\tilde{G}$ is Gorenstein at $0$ and coincides with $P_n$, hence,
$m_0=m_2+\dots+m_{n-1}+u\in \mbox{Int}(n\tilde{G})\cap M$.
Since $m_1+\dots+m_{n-1}$ is in $\mbox{Int}(nQ)$, the point $-e'_n$ is contained in $Q$.
This contradicts to that $m_1$ is a vertex of $Q$.

In the case that $a_n\ge3$, we can decompose $\tilde{G}$ into a union of $a_n-1$
basic $n$-simplices with vertices $\{e'_2, \dots, e'_{n-1}, 
e'_1+(j-1)e'_n, e'_1+je'_n\}$
for $j=1, \dots, a_n-1$.
If  $m_0$ is in the relative interior of the cone of dimension $n-1$
say, $C(e'_2, \dots, e'_{n-1},  e'_1+je'_n)$, then it contradicts by the above reason.
If $m_0$ is contained in the interior of the cone of
one on these $n$-simplices,
say, $C(e'_2, \dots, e'_{n-1}, e'_1+(j-1)e'_n, e'_1+je'_n)$.  Then $m_0=e'_2+
\dots +e'_{n-1} + (e'_1+(j-1)e'_n)+ (e'_1+je'_n)$.  In this case, it contradicts to that
the cone is strictly convex.
Thus the case that $t>n$ does not occur.
\hfill $\Box$

\medskip
\noindent
{\it Proof of Proposition~{\ref{p1a}}}.
We separate our argument into two cases according to the existence of interior lattice
points of $nQ$.

The case I: $\mbox{Int}(nQ)\cap M=\emptyset$.
Then $Q$ is basic and $P\not=Q$.
Let $v'$ be a vertex of $P$ not contained in $Q$.  
As in the proof of Proposition~\ref{p1}, we may set $m_i=e_i$ for $i=1, \dots, n$
and  we can write as 
$$
v'=ae_1+e_2
$$
with a positive integer $a$.
We note that $m_0=e_1+\dots +e_n$.

Set
$$
G:=\mbox{Conv}\{v', e_2, \dots, e_n\}\quad\mbox{and} \quad
\tilde{G}:=\mbox{Conv}\{0, G\}.
$$
Then $m_0$ is contained in the relative interior of $(n-1)G$ and 
in the interior of $n\tilde{G}$.
In other words, $P$ is included in the prism written as 
$x_1\ge0, \dots, x_n\ge0, x_2+\dots+x_n\le1$
by using the coordinates $(x_1, \dots, x_n)$ of $M_{\mathbb{R}}$.

If $v'$ is the only vertex of $P$ other than $m_1, \dots, m_n$, then $a=1$ and
$P$ is isomorphic to $Q_n$ otherwise the vertex $m_3$ is not Gorenstein.

Let $v_j=a_je_1+e_j$ for $j=2, \dots, n$ be vertices of $P$.  
We may set $a_2\ge \dots \ge a_n\ge0$ and $a_3\ge1$.
If $a_{n-1}\ge1$ and $a_n=0$, then the vertex $m_n$ is not Gorenstein.
Thus $P\cong R_n$.

The case II: $ \mbox{Int}(nQ)\not=\emptyset$.
Lemma~\ref{3;l1} shows that $m_0\in \mbox{Int}(nQ)$ and
 that $Q$ coincides with $Q'_n$ (the case (a)) or $P_n$ (the case (b)).
  We have to consider the case that $P\not=Q$.
  
In the case (a) we have 
$$
Q=Q'_n=\mbox{Conv}\{0, e_1, \dots, e_n, e_1+e_2-e_n\}.
$$
We note that $n+1$ vertices of $Q'_n$ except $0$ are on a hyperplane.
Set $\tilde{F}=\mbox{Conv}\{e_1, \dots, e_n, e_1+e_2-e_n\}$.  Since it is not simplex,
$\mbox{Int}(n-1)\tilde{F}\cap M\not=\emptyset$.
If $P\not=Q$, then $\tilde{F}$ is not a facet of $P$.  This contradicts to $\mbox{Int}
(n-1)P\cap M=\emptyset$.   Thus $P=Q$.

In the case (b) we have 
$$
Q=P_n=\mbox{Conv}\{0, 2e_1, e_2, \dots, e_n\}.
$$
We note that $m_0=e_1+\dots +e_n\notin (n-1)P_n$.
Consider the prism $(x_1\ge0, \dots, x_n\ge0, x_2+\dots+x_n\le1)$
by using the coordinates $(x_1, \dots, x_n)$ of $M_{\mathbb{R}}$.
The point $m_0$ is contained in the boundary of the $(n-1)$-tuple of the prism.
We see that if $P\not=Q$, then $P\cong R_n$ as in the case I.
\hfill $\Box$

\section{Fujita's freeness conjecture}\label{sect4}

In this section we give a proof of Theorem~\ref{0;tm1}, which is a 
strong version of Fujino's Theorem\cite{Fn} but restricted to Gorenstein
toric varieties.

We recall the construction of the polarized toric $n$-fold $(X, L)$ from a lattice
$n$-polytope $P$ (see, for instance, \cite{Od} or \cite{Fu}).
For simplicity, we assume that all toric varieties are defined over the complex number field
$\mathbb{C}$.
Let $N$ be a free abelian group of rank $n$ and $M$ the dual with the natural pairing
$\langle, \rangle: M\times N \to \mathbb{Z}$.
Let $T_N:=N\otimes_{\mathbb{Z}}\mathbb{C}^{\times}$ be the algebraic torus of dimension $n$.
Then the group of characters $\mbox{Hom}_{\mbox{gr}}(T_N, \mathbb{C}^{\times})$
can be identified with $M$ and we have $T_N=\mbox{Spec}\ \mathbb{C}[M]$.
Let $P\subset M_{\mathbb{R}}$ be a lattice $n$-polytope.   From $P$ we construct
a polarized toric $n$-fold $(X, L)$ 
satisfying the equality
\begin{equation}\label{4;e1}
\Gamma(X, L) \cong \bigoplus_{m\in P\cap M} \mathbb{C}\ e(m),
\end{equation}
where we write as $e(m)$ the  character corresponding to a lattice point $m\in M$.
Since $X$ contains $T_N$ as an open subset, we can consider $e(m)$ as a rational function on $X$.

For a vertex $v$ of $P$,  let $\sigma(v)\subset N_{\mathbb{R}}$ be the cone dual to
the cone $C_v(P)=\mathbb{R}_{\ge0}(P-v)\subset M_{\mathbb{R}}$. 
Let $\Phi$ be the set of all faces of cones $\sigma(v)$ for all vertices of $P$.
The $\Phi$ is a complete fan in $N$ and defines a toric variety $X$ of dimension $n$.
We note that $X$ is covered by affine open sets $U_v:=\mbox{Spec}\ \mathbb{C}[M\cap
C_v(P)]$.  Here we  define a line bundle $L$ so that
$$
\Gamma(U_v, L) = e(v)\mathbb{C}[M\cap C_v(P)].
$$
Then $L$ is generated by global sections and ample.
By definition $L$ satisfies~(\ref{4;e1}).

 Furthermore, we assume that $X$ is Gorenstein.  For each vertex $v$, then, the cone $C_v(P)$
  contains the lattice point $m_v$ satisfying the equality
  \begin{equation}\label{4;e2}
(\mbox{Int}C_v(P))\cap M = m_v+ C_v(P)\cap M.
\end{equation}
Thus we see that $L+K_X$ is generated by global sections if $P$ contains all $m_v$ in its interior.

\begin{prop}\label{5;p1}
Let $X$ be a projective Gorenstein toric $n$-fold with $n\ge2$.
Let $L$ be an ample line bundle on $X$ satisfying the condition that 
$\Gamma(X, L+K_X)\not=0$ and that the intersection number
$L\cdot C\ge n-1$ for all irreducible invariant curves $C$.
Then $L+K_X$ is nef.
\end{prop}
{\it Proof}.
Let $P\subset M_{\mathbb{R}}$ be the lattice $n$-polytope corresponding to $L$.
The condition $\Gamma(X, L+K_X)\not=0$ is interpreted to $\mbox{Int}P\cap M\not=
\emptyset$.
Since $P$ is Gorenstein, for each vertex $v$ of $P$, 
there exists the lattice point $m_v$ satisfying the equation~(\ref{4;e2}).
From the above observation, it is suffice  to show that the lattice points
$m_v$ are contained in the interior of $P$ for all vertices $v$.

As in the proof of Proposition~\ref{p1}, we take a vertex $v$ of $P$
so that $v$ coincides with the origin $0$ of $M$.
Let $\{m_1, \dots, m_t\}$ be the set of the nearest points on all edges through 
the vertex $v=0$
and $Q:=\mbox{Conv}\{0, m_1, \dots, m_t\}$.
Then $(n-1)Q$ is contained in $P$ because all edges of $P$ have length at least $n-1$.
If $\mbox{Int}(n-1)Q\cap M\not=\emptyset$, then we see that $m_v\in \mbox{Int}(n-1)Q
\subset \mbox{Int}P$ by Lemma~\ref{3;l0}.
We assume that $\mbox{Int}(n-1)Q\cap M=\emptyset$.

If $\mbox{Int}(nQ)\cap M=\emptyset$, then $Q$ is basic, hence, $t=n$ and we may set
$m_i=e_i$ with a $\mathbb{Z}$-basis $\{e_1, \dots, e_n\}$ of $M$.
Then we have $m_v=e_1+\dots +e_n$.
Let $m'$ be a lattice point in the interior of $P$ different from $m_v$.
Then since $\mbox{Conv}\{(n-1)e_1, \dots, (n-1)e_n, m'\}$ contains $m_v$,
the polytope $P$ contains $m_v$ in its interior.

If $\mbox{Int}(nQ)\cap M\not=\emptyset$, then $Q$ is not basic and $m_v\in 
\mbox{Int}(nQ)$ by Lemma~\ref{3;l0}.
From Lemma~\ref{3;l1}, we see that $Q$ is $Q'_n$, that is, $n\ge3$,  $t=n+1$
and $m_i=e_i$ for $i=1, \dots, n$ and $m_{n+1}=e_1+e_2-e_n$ for a $\mathbb{Z}$-basis
$\{e_1, \dots, e_n\}$ of $M$.
In particular, $Q$ has the facet $F$ containing all $m_i$'s and $m_v=e_1+\dots +e_{n-1}$.
Since $m_v\in \mbox{Int}(n-1)F$, the condition $\mbox{Int}P\cap M\not=\emptyset$
implies that $(n-1)Q\not=P$ and $m_v\in \mbox{Int}P$.
\hfill $\Box$
\medskip

By combining Propositions~\ref{p1a} and \ref{5;p1} and Lemma~\ref{l3},
we obtain the proof of Theorem~\ref{0;tm1}.

\section{Fujita's very ampleness conjecture}\label{sect5}

In this section we will give a proof of Theorem~\ref{0;tm2}.
Let $P\subset M_{\mathbb{R}}$ be the Gorenstein 
lattice $n$-polytope corresponding to $L$.
As in the previous section, for a vertex $v$ of $P$ denote $m_v$ the lattice point
satisfying the equation (\ref{4;e1}).
As in the proof of Proposition~\ref{5;p1}, we take a vertex $v$ of $P$
so that $v$ coincides with the origin $0$ of $M$.
Let $\{m_1, \dots, m_t\}$ be the set of the nearest points on all edges through 
the vertex $v=0$
and $Q:=\mbox{Conv}\{0, m_1, \dots, m_t\}$.
We know that $m_v$ is contained in the interior of $(n+1)Q$.

If $\mbox{Int}(r-1)Q\cap M=\emptyset$ and if $\mbox{Int}(rQ)\cap M\not=\emptyset$
for some $r$ with $1\le r\le n$,
then $m_v\in \mbox{Int}(rQ)$ by Lemma~\ref{3;l0} and 
$(n+1-r)Q$ is normal by Lemma~\ref{l1} and 
$$
m_v+(n+1-r)Q\subset \mbox{Int}(n+1)Q \subset \mbox{Int}P.
$$
This implies that $L+K_X$ is very ample on the affine open set $U_v$.

Assume that $\mbox{Int}(n+1)Q\cap M=\emptyset$.  Then $Q$ is basic, $t=n$ and
$m_i=e_i$ for a $\mathbb{Z}$-basis $\{e_1, \dots, e_n\}$ of $M$.
In particular, $Q$ is normal and $m_v=e_1+\dots+e_n$.
Let $F$ be the facet of $Q$ containing all $e_i$.  Then all $m_v+e_i$ are contained in
the relative interior of $(n+1)F$.
Since $P\not=(n+1)Q$ by assumption, all $m_v+e_i$ are contained in the interior of $P$.
This implies that $L+K_X$ is very ample on this $U_v$.
This completes the proof of Theorem~\ref{0;tm2}.

\medskip

We remark that the condition $L\cdot C\ge n+1$ is best possible for all dimension $n\ge2$.
A Gorenstein toric Fano $n$-fold with index $n$ corresponding to $P_n$ or $Q_n$
trivially attains the bound for all $n\ge"$.

Besides Gorenstein toric Fano $n$-fold with index $n$, we also have a Gorenstein toric Fano
$n$-fold with index $n-1$ which attains the bound "$n+1$".

For $n\ge3$, we define  a lattice $n$-simplex  as 
\begin{equation*}
D_n:=\mbox{\rm Conv}\{0, e_1, e_2, e_1+e_2+2e_3, e_4, \dots, e_n\}.
\end{equation*}
Then it is not very ample and satisfies $\mbox{Int}(n-2)D_n\cap M=\emptyset$ and 
$\sharp (\mbox{\rm Int}(n-1)D_n)\cap M=1$, hence, $2D_n$ is normal.
If we denote by $m_0$ the unique interior lattice point of $(n-1)D_n$, then
we have
$$
m_0+(2D_n)\cap M=\mbox{Int}(n+1)D_n\cap M.
$$
Payne has pointed out the case $n=3$ in \cite{P}.

We can characterize $D_n$ among polarized toric varieties as the following way.

\begin{prop}\label{pr1}
Let $D$ be an empty lattice $n$-simplex for $n\ge3$.
Assume that $\sharp (\mbox{\rm Int}(n-1)D)\cap M=1$.
Then $D$ is isomorphic to $D_n$.
\end{prop}
{\it Proof}.
We note that $D$ is not basic.

Consider the case $n=3$.  An empty lattice 3-simplex $D$ is written as
$$
D=\mbox{Conv}\{0, e_1, e_2, e_1+pe_2+qe_3\}
$$
with $1\le p<q$ and $\gcd(p,q)=1$ by \cite{Od}.
Then $\sharp\mbox{Int}(2D)\cap M=q-1$.  Hence $q=2$ and $D\cong D_3$.

Set $n\ge4$.  We may write as
$$
D=\mbox{Conv}\{0, m_1, \dots, m_n\}
$$
for linearly independent $m_1, \dots, m_n\in M$.

First we assume that the facet $F=\mbox{Conv}\{0, m_1, \dots, m_{n-1}\}$ is basic, that is,
$m_i=e_i$ for $i=1, \dots, n-1$.  Then $\mbox{Int}(n-2)F\cap M=\emptyset$.
We may write as
$$
m_n =a_1e_1+\dots +a_ne_n
$$
with $a_i\ge0$ for $i=1, \dots, n-1$ and  $a_n\ge1$.
Set $H$ the hyperplane containing $F$.
Moreover, we assume that all facets of $D$ are basic.
 Then the $n$-th coordinates of lattice points in
$(n-2)D\setminus H$ are at least $a$.
On the other hand, since $\sharp (\mbox{\rm Int}(n-1)D)\cap M=1$,
$(n-2)D$ does not contain lattice points in its interior, hence, $2D$ is normal
by Lemma~\ref{l1}.  Since $2\le n-2$, $(n-2)D$ is also normal, hence,
there exists a lattice point in it whose $n$-th coordinate is 1.  Then $a=1$.
This implies that $D$ is basic and contradicts to the assumption.
Then we see that at least one facet of $D$ is not basic.

Next, if the facet $F$ is not basic, then $(n-2)F$ contains lattice points in its interior,
in fact, interior lattice point is unique point in $\mbox{Int}(n-1)D$.
By induction on $n$, we see that  the facet $F$ is isomorphic to $D_{n-1}$
and $\mathbb{Z}(F\cap M) =\mathbb{Z}e_1+\dots+ \mathbb{Z}e_{n-1}$.
Since $2D$ is normal, we can take  $m_n=e_n$, hence, $D\cong D_n$.
\hfill $\Box$

\medskip

We remark that "$n-1$" in Proposition~\ref{pr1} is essential.
For $n\ge4$, we define a lattice $n$-simplex as
$$
 P:=\mbox{Conv}\{0, e_1, \dots, e_{n-1}, e_1+e_2+e_3 +3e_n\}.
 $$
We note that $P$ is not Gorenstein.
Set $m':=e_1+\dots +e_{n-1} +2e_n$.
Then we see that $\mbox{Int}(n-2)P\cap M =\{m'\}$.

\end{document}